\documentstyle[12pt]{article}
\newtheorem{theorem}{Theorem}[section]
\newtheorem{lemma}[theorem]{Lemma}

\newtheorem{corollary}[theorem]{Corollary}

\newtheorem{theorem A}[theorem]{Theorem A$\!\!$}
\newtheorem{proposition A}[theorem]{Proposition A$\!\!$}
\newtheorem{lemma A}[theorem]{Lemma A$\!\!$}
\newtheorem{corollary A}[theorem]{Corollary A$\!\!$}

\newcommand{\RR}{{\sf R}\hspace{-2.08mm}{\sf R}}

\newcommand{\Pto}{\stackrel{P}{\to}}
\newcommand{\dto}{\stackrel{d}{\to}}
\newcommand{\vto}{\stackrel{v}{\to}}
\newcommand{\wto}{\stackrel{w}{\to}}

\newcommand{\ssfrac}[2]{\mbox{\footnotesize $\frac{#1}{#2}$}}
\newcommand{\half}{\ssfrac{1}{2}}

\newcommand{\sfrown}{\mbox{$\mbox{\scriptsize $\frown$}$}}
\newcommand{\lfrown}
  {\;\mbox{\raisebox{-1.5mm}{$\stackrel{\textstyle <}{\sfrown}$}}\;}
\newcommand{\efrown}
  {\;\mbox{\raisebox{-1.5mm}{$\stackrel{\textstyle =}{\sfrown}$}}\;}

\newcommand{\gfrown}
  {\;\mbox{\raisebox{-1.5mm}{$\stackrel{\textstyle >}{\sfrown}$}}\;}
\newcommand{\qed}{\nolinebreak \hfill $\Box$}

\begin{document}
\begin{center}
{\Large\bf Lebesgue approximation\\[3mm]
of $(2,\beta)$-superprocesses}\\[5mm]
{\large\em By Xin He, Auburn University}
\end{center}
{\small
\begin{quotation} \noindent

{\bf Abstract:} Let $\xi=(\xi_t)$ be a locally finite $(2,\beta)$-superprocess in
$\RR^d$ with $\beta<1$ and $d>2/\beta$. Then for any fixed $t>0$, the random measure
$\xi_t$ can be a.s. approximated by suitably normalized restrictions of
Lebesgue measure to the $\varepsilon$-neighborhoods of ${\rm supp}\,\xi_t$.
This extends the Lebesgue approximation of Dawson-Watanabe superprocesses.
Our proof is based on a truncation of $(\alpha,\beta)$-superprocesses
and uses bounds and asymptotics of hitting probabilities.
\end{quotation}}

{\footnotesize
\begin{description}
\item{\em Author:}\/ Xin He, Department of Mathematics and
Statistics, Auburn University, 221 Parker Hall, Auburn, AL 36849,
USA, {\sf xzh0005@auburn.edu}
\item{\em Running head:}\/ Approximation of $(2,\beta)$-processes.
\item{\em MSC2010 subject classifications:}\/ 60J68
\item{\em Key words and phrases:}\/ Super-Brownian motion, branching mechanism, symmetric stable process, historical cluster, hitting probability, hitting bound, hitting asymptotic, local finiteness, local extinction, neighborhood measure.\\
\end{description}}

\section{\large\bf Introduction}

By an ($\alpha,\beta$)-{\em superprocess} \/(or {\em $(\alpha,\beta)$-process}, for short) in \/$\RR^d$ we mean a vaguely rcll, measure-valued strong Markov process $\xi=(\xi_t)$ in $\RR^d$ satisfying $E_\mu e^{-\xi_tf}=e^{-\mu v_t}$ for suitable functions $f\geq 0$, where $v=(v_t)$ is the unique solution to the {\em evolution equation}\/ $\dot v=\half\Delta_\alpha v-v^{1+\beta}$ with initial condition $v_0=f$. Here $\Delta_\alpha=-(-\Delta)^{\alpha/2}$ is the fractional Laplacian, $\alpha \in (0,2]$ refers to the spatial motion, and $\beta\in(0,1]$ refers to the branching mechanism. When $\alpha=2$ and $\beta=1$ we get the {\em Dawson--Watanabe superprocess}\/ ({\em DW-process} for short), where the spatial motion is standard Brownian motion. General surveys of superprocesses include the excellent monographs and
lecture notes \cite{D93, Dy94, E00, LG99, L11, P02}.

In this paper we consider superprocesses with possibly infinite initial measures. Indeed, by the additivity property of superprocesses, we can construct the $(\alpha,\beta)$-process $\xi$ with any $\sigma$-finite initial measure $\mu$. In Lemma \ref{2.locfinite} we show that $\xi_t$ is a.s.\ locally finite for
every $t>0$ iff $\mu p_\alpha(t,\cdot)<\infty$ for all $t$, where $p_\alpha(t,x)$ denotes the transition density of a symmetric $\alpha$-stable process. Note that when $\alpha=2$, $p_2(t,x)=p_t(x)$ is the standard normal density.

For any measure $\mu$ on $\RR^d$ and constant $\varepsilon>0$, write $\mu^\varepsilon$ for the
restriction of Lebesgue measure $\lambda^d$ to the $\varepsilon
$-neighborhood of ${\rm supp}\,\mu$. For a DW-process $\xi$ in $\RR^d$ with any finite initial measure, Tribe \cite{T94} showed that $\varepsilon ^{2-d}\,\xi^\varepsilon_t\wto c_d\,\xi_t$ a.s.\ as $\varepsilon\to 0$ when $d\geq 3$, where $\wto$ denotes weak convergence and $c_d>0$ is a constant depending on $d$. For a locally finite DW-process $\xi$ in $\RR^2$, Kallenberg \cite{K08} showed that $\tilde m(\varepsilon)\,|\log\varepsilon|\, \xi_t^\varepsilon \vto\xi_t$ a.s.\ as $\varepsilon\to 0$, where $\vto$ denotes vague convergence and $\tilde m$ is a suitable normalizing function. Our main result is Theorem \ref{5.lebesgue}, where we prove that, for a locally finite $(2,\beta)$-process $\xi$ in $\RR^d$ with $\beta<1$ and $d>2/\beta$, $\varepsilon ^{2/\beta-d}\,\xi^ \varepsilon_t\vto c_{\beta,d}\,\xi_t$ a.s.\ as $\varepsilon\to 0$, where $c_{\beta,d}>0$ is a constant depending on $\beta$ and $d$. In particular, the $(2,\beta)$-process $\xi_t$ distributes its mass over ${\rm supp}\,\xi_t$ in a deterministic manner, which extends the corresponding property of DW-processes (cf. \cite{E00}, page 115, or \cite{P02}, page 212). For DW-processes, this property can also be inferred from some deep results involving the exact Hausdorff measure (cf. \cite{DP91}). However, for any $(\alpha,\beta)$-process $\xi$ with $\alpha<2$, ${\rm supp}\,\xi_t=\RR^d$ or $\emptyset$ a.s. (cf. \cite{EP91, P90}), and so the corresponding property fails. Our result shows that this property depends only on the spatial motion.

To prove our main result, we adapt the probabilistic
approach for DW-processes from \cite{K08}. However, the finite variance of DW-processes plays a crucial role there.
In order to deal with the infinite variance of $(2,\beta)$-processes with $\beta<1$, we use a truncation of $(\alpha,\beta)$-processes from \cite{MV07}, which will be further developed in Section 2 of the present paper. By this truncation we may reduce our discussion to the truncated processes, where the variance is finite.

To adapt the probabilistic approach from \cite{K08} to study the truncated processes, we also need to develop some technical tools. Thus, in Section 3 we improve the upper bounds of hitting probabilities for $(2,\beta)$-processes with $\beta<1$ and their truncated processes. As an immediate application, in Theorem \ref{3.extinction} we improve some known extinction criteria of the $(2,\beta)$-process $\xi$ by showing that the local extinction property $\xi_t\dto 0$ and the seemingly stronger support property ${\rm supp}\,\xi_t\dto\emptyset$ are equivalent. Then in Section 4 we derive some asymptotic results of these hitting probabilities. In particular, for the $(2,\beta)$-process $\xi$ we show in Theorem \ref{4.asymhit} that $\varepsilon^{2/\beta-d}P_\mu\{\xi_tB_x^\varepsilon>0\}\to c_{\beta,d}\,(\mu*p_t)(x)$, where $B_x^r$ denotes an open ball around $x$ of radius $r$, which extends the corresponding result for DW-processes (cf. Theorem 3.1(b) in \cite{DIP89}). Since the truncated processes don't have the scaling properties of the $(2,\beta)$-process, our general method is first to study the $(2,\beta)$-process, then to estimate the truncated processes by the $(2,\beta)$-process, in order to get the needed results for the truncated processes.

The extension of results of DW-processes to general $(\alpha,\beta)$-processes is one of the major themes in the research of superprocesses. Since the spatial motion of the $(\alpha,\beta)$-process is not continuous when $\alpha<2$ and the $(\alpha,\beta)$-process has infinite variance when $\beta<1$, many extensions are not straightforward, and some may not even be valid. However, it turns out that several properties of the support of $(2,\beta)$-processes depend only on the spatial motion. These properties include short-time propagation of the support (cf. Theorem 9.3.2.2 in \cite{D93}) and Hausdorff dimension of the support (cf. Theorem 9.3.3.5 in \cite{D93}). Our result also belongs to that category.

In this paper we are mainly using the notations in \cite{K08}. Especially we use relations such as $\efrown$, $\lfrown$, and $\asymp$, where the first two mean
equality and inequality up to a constant factor, and the last one is the combination of $\lfrown$ and $\gfrown$. Other notation will be explained whenever it occurs.

\section{\large\bf Truncated superprocesses and local finiteness}

It is well known that the $(\alpha,1)$-process has weakly continuous sample paths. By contrast, the $(\alpha,\beta$)-process $\xi$ with $\beta<1$ has only weakly rcll sample paths with jumps of the form $\Delta\xi_t=r\delta_x$, for some $t>0$, $r>0$, and $x\in\RR^d$. Let
$$N_{\xi}(dt,dr,dx)=\sum_{(t,r,x):\ \Delta\xi_t=r\delta_x}\delta_{(t,r,x)}.$$
Clearly the point process $N_{\xi}$ on $\RR_+\times\RR_+\times\RR^d$ records all information about the jumps of $\xi$. By the proof of Theorem 6.1.3 in \cite{D93}, we know that $N_{\xi}$ has compensator measure
\begin{equation}
\hat{N}_{\xi}(dt,dr,dx)=c_\beta (dt)r^{-2-\beta}(dr)\xi_t(dx), \label{compensator}
\end{equation}
where $c_\beta$ is a constant depending on $\beta$.
Due to all the ``big" jumps, $\xi_t$ has infinite variance. Some methods for $(\alpha,1)$-processes, which rely on the finite variance of the processes, are not directly applicable to $(\alpha,\beta$)-processes with $\beta<1$.

In \cite{MV07}, Mytnik and Villa introduced a truncation method for ($\alpha,\beta$)-processes with $\beta<1$, which can be used to study ($\alpha,\beta$)-processes with $\beta<1$, especially to extend results of $(\alpha,1)$-processes to ($\alpha,\beta$)-processes with $\beta<1$. Specifically, for the ($\alpha,\beta$)-process $\xi$ with $\beta<1$, we define the stopping time $\tau_{K}=\inf\{t>0:\|\Delta\xi_t\|>K\}$ for any constant $K>0$. Clearly $\tau_{K}$ is the time when $\xi$ has the first jump greater than $K$. For any finite initial measure $\mu$, they proved that one can define $\xi$ and a weakly rcll, measure-valued Markov process $\xi^K$ on a common probability space such that $\xi_t=\xi^K_t$ for $t<\tau_{K}$. Intuitively, $\xi^K$ euqals $\xi$ minus all masses produced by jumps greater than $K$ along with the future evolution of those masses. In this paper, we call $\xi^K$ the {\em truncated K-process} of $\xi$. Since all ``big" jumps are omitted, $\xi^K_t$ has finite variance. They also proved that $\xi^K_t$ and $\xi_t$ agree asymptotically as $K \rightarrow \infty$. We give a different proof of this result, since similar ideas will also be used at several crucial stages later. We write $P_\mu\{\xi\in\cdot\}$ for the distribution of $\xi$ with initial measure $\mu$.

\begin{lemma}\label{2.asymsame}
Fix any finite $\mu$ and $t>0$. Then $P_\mu\{\tau_{K}>t\}\rightarrow 1$ as $K\rightarrow\infty$.
\end{lemma}

{\em Proof:} If $\tau_{K}\leq t$, then $\xi$ has at least one jump greater than $K$ before time $t$. Noting that $N_{\xi}([0,t],(K,\infty),\RR^d)$ is the number of jumps greater than $K$ before time $t$, we get by Theorem 25.22 of \cite{K02} and (\ref {compensator}),
\begin{eqnarray*}
P_\mu\{\tau_{K}\leq t\}&\leq& E_\mu N_{\xi}\!\left([0,t],(K,\infty),\RR^d\right)\\
&=& E_\mu\hat{N}_{\xi}\!\left([0,t],(K,\infty),\RR^d\right) \\
&\efrown& K^{-1-\beta}E_\mu\!\int_0^t\|\xi_s\|ds=t\|\mu\|K^{-1-\beta}\rightarrow 0
\end{eqnarray*}
as $K\rightarrow\infty$, where the last equation holds by $E_\mu\|\xi_s\|=\|\mu\|$.
\qed\\

Using Lemma 1 of \cite{MV07} and a recursive construction, we can prove that $\xi^K_t(\omega)\leq\xi_t(\omega)$ for any $t$ and $\omega$. So indeed, $\xi^K$ is a ``truncation" of $\xi$.

\begin{lemma}\label{2.domination}
We can define $\xi$ and $\xi^K$ on a common probability space such that: \vspace{-1.5mm}
\begin{itemize}
\item[{\em(i)}] $\xi$ is an $(\alpha,\beta)$-process with $\beta<1$ and a finite initial measure $\mu$, and $\xi^K$ is its truncated $K$-process, \vspace{-2mm}
\item[{\em(ii)}] $\xi_t(\omega)\geq\xi^K_t(\omega)$ for any $t$ and $\omega$, \vspace{-2mm}
\item[{\em(iii)}] $\xi_t(\omega)=\xi^K_t(\omega)$ for $t<\tau_K(\omega)$.\vspace{-1mm}
\end{itemize}
\end{lemma}

{\em Proof:}
Let $\xi_{m,n}(t)$ denote the process $\xi_{m,n}$ at time $t$. Use $D([0,\infty),\hat{\mathcal{M}}_d)$ as our $\Omega$, the space of rcll functions from $[0,\infty)$ to $\hat{\mathcal{M}}_d$. We endow $\Omega$ with the Skorohod $J_1$-topology. Let $\mathcal{A}=B$$(\Omega)$.

Let $\zeta_1(t,\omega)=\omega(t)$ be an ($\alpha,\beta$)-process defined on $(\Omega,\mathcal{A},$$P)$ with initial measure $\mu$, and define $\tau_{K_1}=\inf\{t>0:\|\Delta\zeta_1(t)\|>K\}$. Then define a kernel $u$ from $\hat{\mathcal{M}}_d$ to $\Omega$ such that $u(\nu,\cdot)$ is the distribution of an $(\alpha,\beta)$-process with initial measure $\nu$, and a kernel $u^K$ from $\hat{\mathcal{M}}_d$ to $\Omega$ such that $u^K(\nu,\cdot)$ is the distribution of the truncated $K$-process of an $(\alpha,\beta)$-process with initial measure $\nu$. By Lemma 6.9 in \cite{K02}, we can define $\zeta_{1,\infty}$ to be an ($\alpha,\beta$)-process with initial measure $\zeta_1(\tau_{K_1})$ on an extension of $(\Omega,\mathcal{A},$$P)$, and $\zeta_{1,\infty}'$  to be the truncated $K$-process of an $(\alpha,\beta)$-process with initial measure $\zeta_1(\tau_{K_1}^-)$. Now define $\xi_1$ and $\xi_1^K$ by
\[
\xi_1(t) = \left\{
\begin{array}{ll}
\zeta_1(t), & t<\tau_{K_1},\\
\zeta_{1,\infty}(t-\tau_{K_1}), & t\geq\tau_{K_1},
\end{array} \right.
\]
\[
\xi_1^K(t) = \left\{
\begin{array}{ll}
\zeta_1(t), & t<\tau_{K_1},\\
\zeta_{1,\infty}'(t-\tau_{K_1}), & t \geq\tau_{K_1}.
\end{array} \right.
\]
By the strong Markov property of ($\alpha,\beta$)-processes and the above construction, we can verify that $\xi_1$ is an ($\alpha,\beta$)-process. By Lemma 1 in \cite{MV07}, $\xi_1^K$ is the truncated $K$-process of an $(\alpha,\beta)$-process. Moreover, $\xi_1$ and $\xi_1^K$ satisfy conditions (ii) and (iii) on $[0,\tau_{K_1})$.

Let $u'$ be a kernel from $\hat{\mathcal{M}}_d\times\hat{\mathcal{M}}_d$ to $\mathcal{A}\times\mathcal{A}$ such that $u'(\nu,\nu',\cdot,\cdot)$ is the distribution of a pair of two independent $(\alpha,\beta)$-processes with initial measures $\nu$ and $\nu'$ respectively. Define $(\zeta_{2,0},\zeta_{2,1})$ with distribution $$u'\left(\xi^K_1(\tau_{K_1}^-),\xi_1(\tau_{K_1})-\xi^K_1(\tau_{K_1}^-),\cdot,\cdot\right).$$
Let $\zeta_2=\zeta_{2,0}+\zeta_{2,1}$, $\zeta_2'=\zeta_{2,0}$, and $\tau_{K_2}=\inf\{t>0:\|\Delta\zeta_2(t)\|>K\}$. Let $\zeta_{2,\infty}$ be an ($\alpha,\beta$)-process with initial measure $\zeta_2(\tau_{K_2})$, and let $\zeta_{2,\infty}'$ be the truncated $K$-process of an $(\alpha,\beta)$-process with initial measure $\zeta_2'(\tau_{K_2}^-)$.
Now define $\xi_2$ and $\xi_2^K$ by
\[
\xi_2(t) = \left\{
\begin{array}{ll}
\xi_1(t), & t<\tau_{K_1},\\
\zeta_2(t-\tau_{K_1}), & \tau_{K_1}\leq t<\tau_{K_1}+\tau_{K_2},\\
\zeta_{2,\infty}(t-\tau_{K_1}-\tau_{K_2}), & t \geq \tau_{K_1}+\tau_{K_2},
\end{array} \right.
\]
\[
\xi_2^K(t) = \left\{
\begin{array}{ll}
\xi^K_1(t), & t<\tau_{K_1},\\
\zeta_2'(t-\tau_{K_1}), & \tau_{K_1}\leq t<\tau_{K_1}+\tau_{K_2},\\
\zeta_{2,\infty}'(t-\tau_{K_1}-\tau_{K_2}), & t \geq \tau_{K_1}+\tau_{K_2}.
\end{array} \right.
\]
Similarly, $\xi_2$ is an ($\alpha,\beta$)-process and $\xi_2^K$ is the truncated $K$-process of an $(\alpha,\beta)$-process. They satisfy conditions (ii) and (iii) on $[0,\tau_{K_1}+\tau_{K_2})$.

Continue the above construction: For every $n$, define $\xi_n$ and $\xi_n^K$ such that $\xi_n$ is an ($\alpha,\beta$)-process, $\xi_n^K$ it the truncated $K$-process of an ($\alpha,\beta$)-process, and they satisfy conditions (ii) and (iii) on $[0,\sum_{k=1}^n\tau_{K_k})$.


It suffices to prove that $\sum_{k=1}^\infty\tau_{K_k}=\infty$ a.s. Suppose $P(\sum_{k=1}^\infty\tau_{K_k}<\infty)>0$. Then there exist $t$ and $a$ such that $P(\sum_{k=1}^\infty\tau_{K_k}<t)=a>0$. Since for every $n$, $\xi_n$ is an $(\alpha,\beta)$-process with initial measure $\mu$, we get
$$an\leq E_\mu\hat{N}_{\xi_n}\!\left([0,t],(K,\infty),\RR^d\right).$$
Noting that by (\ref {compensator}) $E_\mu\hat{N}_{\xi_n}([0,t],(K,\infty),\RR^d)$ is the same finite constant for different $n$, we get a contradiction. So $\sum_{k=1}^\infty\tau_{K_k}=\infty$ a.s.
\qed\\

Just as the DW-process, the $(\alpha,\beta)$-process $\xi$ and its truncated $K$-process $\xi^K$ also have cluster structures (cf. Section 3 in \cite{DP91}). Specifically, for any fixed $t$, $\xi_t$ is a Cox cluster process, such that the ``ancestors" of $\xi_t$ at time $s=t-h$ form a Cox process directed by $(\beta h)^{-1/\beta}\xi_s$, and the generated $h$-clusters $\eta_h^i$ are conditionally independent and identically distributed apart from shifts. For the truncated $K$-process $\xi^K$, the situation is similar, except that the clusters are different (because of the truncation) and the term $(\beta h)^{-1/\beta}$ for $\xi$ needs to be replaced by $a_K(h)$ (or $a_h$, when $K$ is fixed). Use $\eta_h^{K,i}$ (or $\eta_h^{Ki}$) to denote the generated $h$-clusters of $\xi^K$. Write $P_x\{\eta_t\in\cdot\}$ for the distribution of $\eta_t$ centered at $x\in\RR^d$, and define $P_\mu\{\eta_t\in\cdot\}=\int\mu(dx)P_x\{\eta_t\in\cdot\}$. The following comparison of $a_K(h)$ and $(\beta h)^{-1/\beta}$,  although not used explicitly in the present paper, should be useful in other applications of the truncation method.

\begin{lemma} \label{2.trunances}
Fix any $K>0$. Then as $h\rightarrow 0$,
$$(\beta h)^{1/\beta}\leq a_K(h)\leq 2(\beta h)^{1/\beta}.$$
\end{lemma}

{\em Proof:} From Lemma 3.4 of \cite{DP91} we know that
$$(\beta h)^{1/\beta}=\lim_{\theta\rightarrow\infty}1/v_0(h,\theta),$$
where $v_0(h,\theta)$ is the solution of $\dot{v}=-v^{1+\beta}$ with initial condition $v\equiv \theta$, and
$$A_K(h)=\lim_{\theta\rightarrow\infty}1/v_1(h,\theta),$$
where $v_1(h,\theta)$ is the solution of (1.12) in \cite{MV07} with initial condition $v\equiv \theta$. Define $M_K(\lambda)=C_\beta(K)\lambda+\Phi^K(\lambda)$, where $C_\beta(K)$ and $\Phi^K$ are such as in (1.12) of \cite{MV07}.
Then $M_K$ satisfies
$$\lambda^{1+\beta}\leq M_K(\lambda)\ \mbox{and}\ \lim_{\lambda\rightarrow \infty}\frac{M_K(\lambda)}{\lambda^{1+\beta}}=1.$$
Clearly it is enough to show that $(1/2)v_0(h,\theta)\leq v_1(h,\theta)\leq v_0(h,\theta)$ as $h\rightarrow 0$ and $\theta\rightarrow \infty$. This follows from the above properties of $M_K$.
\qed\\

Unlike the normal densities, we have no explicit expressions for the transition densities of symmetric $\alpha$-stable processes when $\alpha<2$. However, a simple estimate of $p_\alpha(t,x)$ is enough for our needs.

\begin{lemma} \label{2.trandens}
Let $p_\alpha(t,x)$, $\alpha\in(0,2]$, $t>0$, and $x\in\RR^d$, denote the transition densities of a symmetric $\alpha$-stable process on
$\RR^d$. Then for any fixed $\alpha$ and $d$,
$$p_\alpha(t,x+y)\lfrown p_\alpha(2t,x),\quad |y|^\alpha\leq t.$$
\end{lemma}

{\em Proof:} First let $\alpha=2$. Note that $p_2(t,x)=p_t(x)$ is the standard normal density on $\RR^d$. For $|x|\leq 4\sqrt{t}$, trivially $p_t(x+y)\lfrown p_{2t}(x)$. For $|x|>4\sqrt{t}$, it suffices to check that
$$-\frac{|x+y|^2}{2t}\leq-\frac{|x|^2}{4t},$$
that is, $2|x+y|^2\geq |x|^2$, which follows easily from $|x|\geq 4|y|$.

Now let $\alpha<2$. By the arguments after Remark 5.3 of \cite{BL02},
\begin{equation}
p_\alpha(t,x)\asymp\left(t^{-d/\alpha}\wedge\frac{t}{|x|^{d+\alpha}}\right). \label{trandens}
\end{equation}
Choose $K>2^{1/\alpha}$ to satisfy
$1\leq 2(1-1/K)^{d+\alpha}$. Since $|y|\leq t^{1/\alpha}$, we have for $|x|> Kt^{1/\alpha}$,
$$\frac{t}{|x+y|^{d+\alpha}}\leq\frac{2t}{|x|^{d+\alpha}}.$$
Noticing also that $(2t)/|x|^{d+\alpha}< (2t)^{-d/\alpha}$
for $|x|> Kt^{1/\alpha}$, we get $p_\alpha(t,x+y)\lfrown p_\alpha(2t,x)$ for $|y|\leq t^{1/\alpha}$ and $|x|> Kt^{1/\alpha}$. The
same inequality holds trivially for $|y|\leq t^{1/\alpha}$ and $|x|\leq Kt^{1/\alpha}$.
\qed\\

Using Lemma \ref{2.domination} and Lemma \ref{2.trandens}, we can generalize Lemma 3.2 in \cite{K08} to any $(\alpha,\beta)$-process.

\begin{lemma} \label{2.locfinite} 
Let $\xi$ be an $(\alpha,\beta)$-process in $\RR^d$, $\alpha\in(0,2]$ and $\beta\in(0,1]$, and fix any $\sigma$-finite measure $\mu$. Then for any
fixed $t>0$, the following two conditions are equivalent:
\vspace{-1.5mm}
\begin{itemize}
\item[{\em(i)}] $\xi_t$ is locally finite a.s.\ $P_\mu$,
\vspace{-2mm}
\item[{\em(ii)}] $E_\mu\xi_t$ is locally finite.
\vspace{-1.5mm}
\end{itemize}
Furthermore, {\em(i)} and\/ {\em(ii)} hold for every $t>0$ iff
\vspace{-1.5mm}
\begin{itemize}
\item[{\em(iii)}] $\mu p_t<\infty$ for all $t>0$,
\vspace{-1.5mm}
\end{itemize}
and if $\alpha<2$, then {\em(iii)} is equivalent to
\vspace{-1.5mm}
\begin{itemize}
\item[{\em(iv)}] $\mu p_t<\infty$ for some $t>0$.
\vspace{-1mm}
\end{itemize}
\end{lemma}

{\em Proof:}\/ The formulas for $E_\mu \xi_t$ and $E_\mu \xi_t^2$ (when $\beta<1$), well known for finite $\mu$, as well as the formulas in Lemma 3 of \cite{MV07}, extend by monotone convergence to any $\sigma$-finite measure $\mu$. We also need the simple inequality that for any fixed $\alpha<2$, $s$, and $t$,
\begin{equation}
p_\alpha(s,x)\asymp p_\alpha(t,x). \label{simple}
\end{equation}
To prove it, use (\ref {trandens}) and consider three cases: $|x|\leq (s\wedge t)^{1/\alpha}$, $|x|\geq (s\vee t)^{1/\alpha}$, and $(s\wedge t)^{1/\alpha}<|x|<(s\vee t)^{1/\alpha}$.

If $\alpha=2$ and $\beta=1$, then this is Lemma 3.2 of \cite{K08}. For $\alpha<2$ and $\beta=1$, using Lemma \ref{2.trandens} and (\ref {simple}) we can proceed as in Lemma 3.2 of \cite{K08}. For example, for any fixed $t>0$ and $x\in\RR^d$, $p(t,x-u)\lfrown p(|x|^{1/\alpha},x-u)\lfrown p(2|x|^{1/\alpha},-u)=p(2|x|^{1/\alpha},u)$ yields $\mu*p(t,\cdot)(x)<\infty$.

Now assume $\beta<1$. Condition (ii) clearly implies (i). Conversely, suppose that $E_\mu\xi_t B=\infty$ for some $B$. Then $E_\mu\xi^K_t B=\infty$ for any fixed $K>0$ by Lemma 3 of \cite{MV07}. Also, we get by Lemma 3 of \cite{MV07},
$$
P_\mu\!\left\{\frac{\xi^K_tB}{E_\mu\xi^K_tB}>r\right\}
 \geq(1-r)^2\,\frac{(E_\mu\xi^K_tB)^2}{E_\mu(\xi^K_tB)^2}
 \geq\frac{(1-r)^2}{1+ct\,(E_\mu\xi^K_tB)^{-1}}\,
$$
for any $r\in(0,1)$.
Arguing as in the proof of Lemma 3.2 in \cite{K08}, we get $\xi^K_t B=\infty$ a.s., and so $\xi_t B=\infty$ a.s. by Lemma \ref{2.domination}. In particular, this shows that (i) implies (ii). To prove the equivalence of (ii) and (iii), again using Lemma \ref{2.trandens} and (\ref {simple}) we can proceed as in Lemma 3.2 of \cite{K08}. The last assertion is obvious from (\ref {simple}).
\qed\\

\section{\large\bf Hitting bounds and neighborhood measures}

The Lebesgue approximation depends crucially on estimates of the hitting probability $P_\mu\{\xi_t B_0^\varepsilon>0\}$. In this section, we first estimate $P_\mu\{\xi_t B_0^\varepsilon>0\}$ and $P_\mu\{\xi^K_t B_0^\varepsilon>0\}$. Then we use these estimates to study multiple hitting and neighborhood measures of the clusters $\eta^K_h$ associated with the truncated $K$-process $\xi^K$. We begin with a well-known relationship between the hitting probabilities of $\xi_t$ and $\eta_t$, which can be proved as in Lemma 4.1 of \cite{K08}.

\begin{lemma} \label{3.hitcomp}
Let the $(\alpha,\beta)$-process $\xi$ in $\RR^d$ with associated clusters
$\eta_t$ be locally finite under $P_\mu$, let $\xi^K$ be its truncated $K$-process with associated clusters
$\eta^K_t$, and fix any $B\in{\cal B}^d$. Then
\begin{eqnarray*}
P_\mu\{\eta_tB>0\}&=&-\,(\beta t)^{1/\beta}\,\log\,(1-P_\mu\{\xi_tB>0\}),\\
P_\mu\{\xi_tB>0\}&=& 1-\exp\,\left(-(\beta t)^{-1/\beta}P_\mu\{\eta_tB>0\}\right),\\
P_\mu\{\eta^K_tB>0\}&=&-\,a_t\,\log\,(1-P_\mu\{\xi^K_tB>0\}),\\
P_\mu\{\xi^K_tB>0\}&=& 1-\exp\,(-a^{-1}_t P_\mu\{\eta^K_tB>0\}).
\end{eqnarray*}
In particular, $P_\mu\{\xi_tB>0\}\sim (\beta t)^{-1/\beta}P_\mu\{\eta_tB>0\}$ and $P_\mu\{\xi^K_tB>0\}\sim a^{-1}_t P_\mu\{\eta^K_tB>0\}$ as
either side tends to 0.

\end{lemma}

Upper and lower bounds of $P_\mu\{\xi_t B_0^\varepsilon>0\}$ have been obtained by Delmas \cite{De99}, using the Brownian snake. However, in this paper we need the following improved upper bound.

\begin{lemma} \label{3.hitprob}
Let $\eta_t$ be the clusters of a $(2,\beta)$-process $\xi$ in $\RR^d$ with $\beta<1$ and $d>2/\beta$, let $\eta^K_t$ be the clusters of $\xi^K$, the truncated $K$-process of $\xi$, and consider a $\sigma$-finite measure $\mu$ on $\RR^d$. Then for $0<\varepsilon\leq\sqrt t$,
\begin{itemize}
\item[{\em(i)}]
$\displaystyle
\mu p_{t'}
 \lfrown \varepsilon^{2/\beta-d}(\beta t)^{-1/\beta}P_\mu\{\eta_tB_0^\varepsilon>0\}
 \lfrown\mu p_{2t},
$ where $t'=\beta t/(1+\beta)$,

\item[{\em(ii)}]
$\displaystyle
\varepsilon^{2/\beta-d}a^{-1}_tP_\mu\{\eta^K_tB_0^\varepsilon>0\}
 \lfrown\mu p_{2t}.
$

\end{itemize}
\end{lemma}

{\em Proof:} (i) From the proof of Theorem 2.3 in \cite{De99} we know that
$$
P_x\{\xi_t B_0^\varepsilon>0\}=1-\exp(-N_x\{Y_tB_0^\varepsilon>0\}),
$$
where $N_x$ and $Y_t$ are defined in Section 4.2 of \cite{De99}.
Comparing this with Lemma \ref{3.hitcomp} yields
$$
(\beta t)^{-1/\beta} P_x\{\eta_tB_0^\varepsilon>0\}=N_x\{Y_tB_0^\varepsilon>0\}.
$$
By Proposition 6.2 in \cite{De99} we get the lower bound. For our upper bound, we will now improve the upper bound in Proposition 6.1 of \cite{De99}.

For $0<\varepsilon/2<\sqrt{t}$, define
\begin{eqnarray*}
\Delta&=&\{(r,y)\in\RR^+\times\RR^d,\ r<t, |y|>\varepsilon/2\}\\
&&\bigcup\{(r,y)\in\RR^+\times\RR^d,\ r<t-\varepsilon^2/16, |y|\leq \varepsilon/2\}.
\end{eqnarray*}
Following the proof of Proposition 6.1 in \cite{De99}, we have
$$(\beta t)^{-1/\beta} P_x\{\eta_tB_0^\varepsilon>0\}\lfrown \varepsilon^{-2/\beta}P_0\{\gamma_s\in B_x^{\varepsilon/2}\ \mbox{for some } s\in [t-\varepsilon^2/16,t)\},$$
where $\gamma$ is a standard Brownian motion. Define
$$T=\inf \{s\geq t-\varepsilon^2/16: \gamma_s\in B_x^{\varepsilon/2}\}.$$
Then $\{T<t\}=\{\gamma_s\in B_x^{\varepsilon/2}\ \mbox{for some } s\in [t-\varepsilon^2/16,t)\}$. To get our upper bound, it remains to show that
$$
P_0\{T<t\}\lfrown\varepsilon^d p_{2t}(x).
$$

To prove this, we need the elementary fact that for any $x\in\RR^d$, $\varepsilon>0$, $y\in B_x^{\varepsilon/2}$, and $s\leq s'=\varepsilon^2/16$,
$$
P_y\{\gamma_s\notin B_x^{\varepsilon}\}\leq P_z\{\gamma_{s'}\notin B_x^{\varepsilon}\}\lfrown P_z\{\gamma_{s'}\in B_x^{\varepsilon}\}\leq P_y\{\gamma_s\in B_x^{\varepsilon}\},
$$
where $z$ is a point on the surface of $B_x^{\varepsilon/2}$, and the third relation holds since $P_z\{\gamma_{s'}\notin B_x^{\varepsilon}\}$ and
$P_z\{\gamma_{s'}\in B_x^{\varepsilon}\}$ are both positive constants. Now return to $P_0\{T<t\}$. Noting $t-T\leq\varepsilon^2/16$ on $\{T<t\}$, we get
\begin{eqnarray*}
P_0\{T<t\}&=& P_0\{T<t,\gamma_t\in B_x^{\varepsilon}\}+P_0\{T<t,\gamma_t\notin B_x^{\varepsilon}\}\\
&=&P_0\{T<t,\gamma_t\in B_x^{\varepsilon}\}+P_0\{P_{\gamma_T}\{\gamma_{t-T}\notin B_x^{\varepsilon}\},T<t\}\\
&\lfrown& P_0\{T<t,\gamma_t\in B_x^{\varepsilon}\}+P_0\{P_{\gamma_T}\{\gamma_{t-T}\in B_x^{\varepsilon}\},T<t\}\\
&=& P_0\{T<t,\gamma_t\in B_x^{\varepsilon}\}+P_0\{T<t,\gamma_t\in B_x^{\varepsilon}\}\\
&\lfrown& P_0\{\gamma_t\in B_x^{\varepsilon}\}\lfrown \varepsilon^d p_{2t}(x),
\end{eqnarray*}
where the second and fourth relations hold by the strong Markov property of Brownian motion and the last relation holds by Lemma \ref{2.trandens}.

(ii) This is obvious from (i), Lemma \ref{2.domination}, and Lemma \ref{3.hitcomp}.
\qed\\

As an immediate application of the improved upper bound, we may improve some known extinction criteria for $(2,\beta)$-processes in $\RR^d$ with $\beta<1$ and $d>2/\beta$. This extends Theorem 4.5 of \cite{K08} for DW-processes of dimension $d\geq 2$. Note that the properties $\xi_t\dto 0$ and ${\rm supp}\,\xi_t\dto\emptyset$ are given by $\xi_t B\Pto 0$ and $1\{\xi_t B>0\}\Pto 0$, respectively, for any bounded Borel set $B$.

\begin{theorem} \label{3.extinction}
Let $\xi$ be a locally finite $(2,\beta)$-process in $\RR^d$, $\beta<1$ and $d>2/\beta$,
with arbitrary initial distribution. Then these conditions are
equivalent \vspace{-2mm} as $t\to\infty$:
\begin{itemize}
\item[{\em(i)}] $\xi_t\dto 0$,
 \vspace{-3mm}
\item[{\em(ii)}] ${\rm supp}\,\xi_t\dto\emptyset$,
 \vspace{-2mm}
\item[{\em(iii)}] $\xi_0p_t\Pto 0$.
\end{itemize}
\end{theorem}

{\em Proof:} By Lemma \ref{3.hitcomp} and Lemma \ref{3.hitprob}(i) we get for any fixed $r$
$$
P_\mu\{\xi_tB_0^r>0\}\leq (\beta t)^{-1/\beta}P_\mu\{\eta_tB_0^r>0\}\lfrown \mu p_{2t},
$$
and so $P_\mu\{\xi_tB_0^r>0\}\lfrown \mu p_{2t}\wedge 1$.
For a general initial distribution,
$$P\{\xi_tB_0^r>0\}\lfrown E(\xi_0 p_{2t}\wedge 1),$$
which shows that (iii) implies (ii). Since clearly (ii) implies (i), it remains to prove that (i) implies (iii).

Let $\xi$ be locally finite under
$P_\mu$. We first choose $f\in C_c^{++}(\RR^d)$ with ${\rm supp}f\in B_0^1$, where $C_c^{++}(\RR^d)$ is such as in Proposition 2.6 of \cite{Kl98}. Clearly $\xi_t f\Pto 0$ if $\xi_t B_0^1\Pto 0$. By dominated convergence
$$\exp (-\mu v_t)=E_\mu \exp(-\xi_t f)\rightarrow 1,$$
and so $\mu v_t\rightarrow 0$. By Proposition 2.6 of \cite{Kl98}, we have for $t$ large enough
$$
p_{t/2}(x) \asymp \phi(t/2,x) \lfrown v_t(x),
$$
and so $\mu p_{t/2}\rightarrow 0$. For general $\xi_0$, we may proceed as in the proof of Theorem 4.5 in \cite{K08}.
\qed\\

The following simple fact is often useful to extend results for
finite initial measures $\mu$ to the general case.

\begin{lemma} \label{3.initial}
Let the $(2,\beta)$-process $\xi$ in $\RR^d$ with $\beta<1$ and $d>2/\beta$ be locally finite under
$P_\mu$, and suppose that $\mu\geq\mu_n\downarrow 0$. Then
$P_{\mu_n}\{\xi_tB>0\}\to 0$ as $n\to\infty$ for any fixed $t>0$ and
$B\in \hat{\cal B}^d$.
\end{lemma}

{\em Proof:}
Follow the proof of Lemma 4.3 in \cite{K08}, then use Lemma \ref{2.locfinite}, Lemma \ref{3.hitcomp}, and Lemma \ref{3.hitprob}(i).
\qed\\

As in \cite{K08} we need to estimate the probability that a ball in $\RR^d$ is hit by more than one subcluster of the truncated $K$-process $\xi^K$. This is where the truncation of $\xi$ is needed.

\begin{lemma} \label{3.multhit}
Fix any $K>0$. Let $\xi^K$ be the truncated $K$-process of a $(2,\beta)$-process $\xi$ in $\RR^d$ with $\beta<1$ and $d>2/\beta$.
For any $t\geq h>0$ and $\varepsilon>0$, let $\kappa_h^\varepsilon$ be the number of $h$-clusters of $\xi^K_t$ hitting
$B_0^\varepsilon$ at time $t$. Then for $\varepsilon^2\leq h\leq t$,
$$
E_\mu\kappa_h^\varepsilon(\kappa_h^\varepsilon-1)\lfrown
 \varepsilon^{2(d-2/\beta)}\left(\,h^{1-d/2}\mu p_t+(\mu
 p_{2t})^2\,\right).
$$
\end{lemma}

{\em Proof:}
Follow Lemma 4.4 in \cite{K08}, then use Lemma 3 of \cite{MV07} and Lemma \ref{3.hitprob}(ii).
\qed\\

Now we consider the neighborhood measures of the clusters $\eta^K_h$ associated with the truncated $K$-process $\xi^K$.
For any measure $\mu$ on $\RR^d$ and constant $\varepsilon>0$, we define the associated {\em neighborhood measure} $\mu^\varepsilon$
as the restriction of Lebesgue measure $\lambda^d$ to the $\varepsilon$-neigh\-bor\-hood of ${\rm supp}\,\mu$, so that
$\mu^\varepsilon$ has Lebesgue density \mbox{$1\{\mu B_x^\varepsilon>0\}$}. Let $p_h^{K,\varepsilon}(x)=P_x \{\eta^K_hB_0^\varepsilon>0\}$, where the $\eta^K_h$ are clusters of $\xi^K$. Write $p_h^{K,\varepsilon}(x)=p_h^{K\varepsilon}(x)$ and $(\eta^{K,i}_h)^\varepsilon=\eta^{Ki\varepsilon}_h$ for convenience.

\begin{lemma} \label{3.variance}

Let $\xi^K$ be the truncated $K$-process of a $(2,\beta)$-process $\xi$ in $\RR^d$ with $\beta<1$ and $d>2/\beta$. Let the $\eta_h^{Ki}$ be conditionally independent $h$-clusters of $\xi^K$, rooted at the points of a Poisson process $\zeta$ with $E\zeta=\mu$. Fix any measurable function $f\geq 0$ on $\RR^d$. Then, \vspace{-1mm}
\begin{itemize}
\item[{\em(i)}]
$E_\mu\,{\sum}_i\,\eta_h^{Ki\varepsilon}
 =(\mu*p_h^{K\varepsilon})\cdot\lambda^d$ , \vspace{-1mm}
\item[{\em(ii)}] $E_\mu{\rm Var}\left[{\sum}_i \eta_h^{Ki\varepsilon}f|\zeta\right] \lfrown
 a_h\varepsilon^{d-2/\beta} h^{d/2}\,\|f\|^2\|\mu\|$ for $\varepsilon^2\leq h$.
\end{itemize}
\end{lemma}

{\em Proof:}
(i)
Follow the proof of Lemma 6.2 (i) in \cite{K08}.

(ii)
First,
$$
{\rm Var}_x(\eta_h^{K\varepsilon} f)\leq E_x(\eta_h^{K\varepsilon} f)^2\leq E_x \|\eta_h^{K\varepsilon}\|^2\,\|f\|^2=\|f\|^2 E_x \|\eta_h^{K\varepsilon}\|^2.
$$
For $E_x \|\eta_h^{K\varepsilon}\|^2$, using Cauchy inequality and Lemma \ref{3.hitprob}(ii), we get
\begin{eqnarray*}
E_x \|\eta_h^{K\varepsilon}\|^2
&=& E_x\left(\int 1\{\eta^K_h B^\varepsilon_y>0\}dy\int 1\{\eta^K_h B^\varepsilon_z>0\}dz\right)\\
&=& \int\int P_x\left(\{\eta^K_h B^\varepsilon_y>0\}\cap\{\eta^K_h B^\varepsilon_z > 0\}\right)dydz\\
&\leq& \int\int (P_x\{\eta^K_h B^\varepsilon_y > 0\}P_x\{\eta^K_h B^\varepsilon_z > 0\})^{1/2}dydz\\
&\lfrown& a_h\varepsilon^{d-2/\beta}\int\int (p_{2h}(y-x)p_{2h}(z-x))^{1/2} dydz\\
&\efrown& a_h\varepsilon^{d-2/\beta}h^{d/2}\int\int p_{4h}(y-x)p_{4h}(z-x) dydz\\
&=& a_h\varepsilon^{d-2/\beta} h^{d/2}. \vspace{-3mm}
\end{eqnarray*}
Hence, by independence
$$
E_\mu{\rm Var}\left[{\sum}_i \eta_h^{Ki\varepsilon}f|\zeta\right]=E_\mu\int\zeta(dx){\rm Var}_x(\eta_h^{K\varepsilon} f)\lfrown a_h\varepsilon^{d-2/\beta}h^{d/2}\,\|f\|^2\|\mu\|.
$$
\qed\\

We also need to estimate the overlap between subclusters.

\begin{lemma} \label{3.overlap}
Let $\xi^K$ be the truncated $K$-process of a $(2,\beta)$-process $\xi$ in $\RR^d$ with $\beta<1$ and $d>2/\beta$. For any fixed $t>0$, let $\eta_h^{Ki}$ denote the subclusters in $\xi^K$ of age $h>0$. Fix any $\mu\in\hat{\cal M}_d$. Then as $\varepsilon^2\leq h\to 0$,
$$
E_\mu\!\left\|\,{\sum}_i\eta_h^{Ki\varepsilon}-\xi_t^{K\varepsilon}\right\|
 \lfrown \varepsilon^{2(d-2/\beta)}h^{1-d/2}.
$$
\end{lemma}
{\em Proof:}
Follow the proof of Lemma 6.3(i) in \cite{K08}, then use Lemma \ref{3.hitprob}(ii).
\qed\\

\section{\large\bf Hitting asymptotics}

For a DW-process $\xi$ of dimension $d\geq 3$, we know from Theorem 3.1(b) of Dawson, Iscoe, and Perkins \cite{DIP89} that, as $\varepsilon\to 0$,
$$
\varepsilon^{2-d}P_\mu\{\xi_tB_x^\varepsilon>0\}\to c_d\,(\mu*p_t)(x),
$$
uniformly for bounded $\|\mu\|$, bounded $t^{-1}$, and $x\in\RR^d$. A similar result for DW-processes of dimension $d=2$ is Theorem 5.3(ii) of \cite{K08}.
 In this section, using Lemma \ref{3.hitprob}(i), we can prove the corresponding result for $(2,\beta)$-processes in $\RR^d$ with $\beta<1$ and $d>2/\beta$.

First we fix a continuous function $f$ on $\RR^d$ such that $0<f(x)\leq1$ for $x\in B_0^1$ and $f(x)=0$ otherwise. Let $v_\lambda$ be the solution of $\dot v=\half\Delta v-v^{1+\beta}$ with initial condition $v(0)=\lambda f$. Since $v_\lambda$ is increasing in $\lambda$, we can define $v_\infty=\lim_{\lambda\rightarrow\infty} v_\lambda$. Using Lemma \ref{3.hitprob}(i), we can get an upper bound of $v_\infty$, similar to Lemma 3.2 in \cite{DIP89}.

\begin{lemma} \label{4.upperbound}
For any $t\geq1$ and $x\in\RR^d$, $v_\infty(t,x)\lfrown p(2t,x)$.

\end{lemma}

{\em Proof:} Letting $\lambda\rightarrow\infty$ in $E_x\exp (-\xi_t \lambda f)=\exp[-v_\lambda(t,x)]$,
we get
$$P_x\{\xi_tB_0^1>0\}=1-\exp[-v_\infty(t,x)].$$
Comparing this with Lemma \ref{3.hitcomp} yields
\begin{equation} \label{vinfinity}
v_\infty(t,x)=(\beta t)^{-1/\beta}P_x\{\eta_tB_0^1>0\}.
\end{equation}
Now Lemma \ref{4.upperbound} follows from Lemma \ref{3.hitprob}(i).
\qed\\

As in Lemma 3.3 of \cite{DIP89}, we can apply a PDE result to get the uniform convergence of $v_\infty$. Notice that the improved upper bound in Lemma \ref{3.hitprob}(i) is crucial here.

\begin{lemma} \label{4.asymsolu}
There exists a constant $c_{\beta,d}>0$ such that
$$
\lim_{\varepsilon\rightarrow 0}\varepsilon^{-d}v_\infty(\varepsilon^{-2}t,\varepsilon^{-1}x)=c_{\beta,d}\cdot p(t,x).
$$
The convergence is uniform for bounded $t^{-1}$ and $x\in\RR^d$.
\end{lemma}

{\em Proof:}
We follow the proof of Lemma 3.3 in \cite{DIP89}. By Lemma \ref{4.upperbound}, $v_\infty(t,x)$ is finite for any $t\geq 1$ and $x\in\RR^d$. Then by a standard regularity argument in PDE theory,
\begin{equation}\label{equation}
\dot v_\infty=\half\Delta v_\infty-v_\infty^{1+\beta}
\end{equation}
on $[1,\infty)\times\RR^d$.
By Lemma \ref{4.upperbound}, $v_\infty(1)\in L^1(\RR^d)$. Set
$$w_\varepsilon(t,x)=\varepsilon^{-d}u_\infty(1+\varepsilon^{-2}t,\varepsilon^{-1}x).$$
Then by (\ref{equation}), $\dot{w}_\varepsilon=\half\Delta w_\varepsilon-\varepsilon^{d-2}w_\varepsilon^{1+\beta}$
with initial condition $w_\varepsilon(0,x)=\varepsilon^{-d}u_\infty(1,\varepsilon^{-1}x)$.

Applying Proposition 3.1 in \cite{GV84} gives
$$\lim_{\varepsilon\rightarrow 0}\varepsilon^{-d}v_\infty(1+\varepsilon^{-2}t,\varepsilon^{-1}x)=c_{\beta,d}\cdot p(t,x),$$
uniformly on compact subsets of $(0,\infty)\times \RR^d$.
Together with Lemma \ref{4.upperbound} this yields the uniform convergence on $[a,\infty)\times \RR^d$ for any $a>0$. Moreover, letting $t=t'-\varepsilon^2$,
we get
$$\lim_{\varepsilon\rightarrow 0}\varepsilon^{-d}v_\infty(\varepsilon^{-2}t',\varepsilon^{-1}x)=c_{\beta,d}\cdot p(t',x),$$
uniformly on $[a,\infty)\times \RR^d$ for any $a>0$.

It remains to prove that $c_{\beta,d}>0$. Using (\ref{vinfinity}) and the lower bound in Lemma \ref{3.hitprob}(i), we obtain
\begin{eqnarray*}
\varepsilon^{-d}v_\infty(\varepsilon^{-2}t,\varepsilon^{-1}x)
&=& \varepsilon^{-d}(\beta t)^{-1/\beta}P_{\varepsilon^{-1}x}\{\eta_{\varepsilon^{-2}t} B_0^1>0\}\\
&\gfrown& \varepsilon^{-d}p\left(\frac{\beta\varepsilon^{-2}t}{1+\beta},\varepsilon^{-1}x\right)=p\left(\frac{\beta t}{1+\beta},x\right),
\end{eqnarray*}
and so $c_{\beta,d}>0$.
\qed\\

Now we can derive the asymptotic hitting rate for a $(2,\beta)$-process.

\begin{theorem} \label{4.asymhit}
Let the $(2,\beta)$-process $\xi$ in $\RR^d$ with $\beta<1$ and $d>2/\beta$ be locally finite under
$P_\mu$. Fix any $t>0$ and $x\in\RR^d$. Then as $\varepsilon\rightarrow 0$,
$$
\varepsilon^{2/\beta-d}P_\mu\{\xi_tB_x^\varepsilon>0\}\to c_{\beta,d}(\mu*p_t)(x).
$$
The convergence is uniform for bounded $\|\mu\|$, bounded $t^{-1}$, and $x\in\RR^d$. Similar results hold for the clusters $\eta_t$ with $p_t$ replaced by $(\beta t)^{1/\beta}p_t$.
\end{theorem}

{\em Proof:} We first prove that as $\varepsilon\rightarrow 0$,
\begin{equation} \label{realone}
\varepsilon^{2/\beta-d}(\beta t)^{-1/\beta}P_\mu\{\eta_tB_x^\varepsilon>0\}\to c_{\beta,d}(\mu*p_t)(x),
\end{equation}
uniformly for bounded $\|\mu\|$, bounded $t^{-1}$, and $x\in\RR^d$.

Use $\mu-x$ to denote the measure $\mu$ shifted by $-x$. If $\mu$ is finite, then by the scaling of $\eta$, (\ref{vinfinity}), and Lemma \ref{4.asymsolu}, we can get the following chain of relations, which proves the uniform convergence of (\ref{realone}):
\begin{eqnarray*}
\lefteqn{\varepsilon^{2/\beta-d}(\beta t)^{-1/\beta}P_\mu\{\eta_tB_x^\varepsilon>0\}}\\
&=&\varepsilon^{2/\beta-d}(\beta t)^{-1/\beta}P_{\mu-x}\{\eta_tB_0^\varepsilon>0\}\\
&=&\varepsilon^{2/\beta-d}(\beta t)^{-1/\beta}\int P_y\{\eta_tB_0^\varepsilon>0\}(\mu-x)(dy)\\
&=&\varepsilon^{2/\beta-d}(\beta t)^{-1/\beta}\int P_{y/\varepsilon}\{\eta_{t/\varepsilon^2}B_0^1>0\}(\mu-x)(dy)\\
&=&\varepsilon^{-d}\int v_\infty(\varepsilon^{-2}t,\varepsilon^{-1}y)(\mu-x)(dy)\rightarrow c_{\beta,d}(\mu*p_t)(x).
\end{eqnarray*}

Let $\mu$ be an infinite $\sigma$-finite measure satisfying $\mu p_t<\infty$ for all $t$. From the proof of Lemma \ref{2.locfinite}, we know that $(\mu*p_{2t})(x)<\infty$ for any $x\in\RR^d$. Then by dominated convergence based on Lemma \ref{3.hitprob}(i), we can still get (\ref{realone}).

Now we turn to $\xi_t$. First note that by Lemma \ref{3.hitcomp}, as $\varepsilon\rightarrow0$,
\begin{eqnarray}
\varepsilon^{2/\beta-d}P_\mu\{\xi_tB>0\}\rightarrow c &\Leftrightarrow& \varepsilon^{2/\beta-d}(\beta t)^{-1/\beta}P_\mu\{\eta_tB>0\}\rightarrow c,\label{cluster}\\
\varepsilon^{2/\beta-d}P_\mu\{\xi^K_tB>0\}\rightarrow c &\Leftrightarrow& \varepsilon^{2/\beta-d}a_t^{-1}P_\mu\{\eta^K_tB>0\}\rightarrow c.\label{kcluster}
\end{eqnarray}
It remains to prove the uniform convergence for $\xi_t$. Since $(\mu*p_t)(x)\leq t^{-d/2}\|\mu\|$,
we know that by (\ref{realone}), $(\beta t)^{-1/\beta}P_\mu\{\eta_tB_x^\varepsilon>0\}\rightarrow 0$,
uniformly for bounded $\|\mu\|$, bounded $t^{-1}$, and $x\in\RR^d$. Then we may use Lemma \ref{3.hitcomp} to get the uniform convergence for $\xi_t$.
\qed\\

The following result, especially part (ii), will play a crucial role in Section 5. Here we approximate the hitting probabilities $p^{K\varepsilon}_h$ by suitably normalized Dirac functions. This will be used in Lemma \ref{5.trunleb} to prove the Lebesgue approximation of $\xi^K$.

\begin{lemma} \label{4.shorttime}
Let $p_h^\varepsilon(x)=P_x \{\eta_hB_0^\varepsilon>0\}$, where
the $\eta_h$ are clusters of a $(2,\beta)$-process $\xi$ in $\RR^d$ with $\beta<1$ and $d>2/\beta$.
Recall that $p_h^{K\varepsilon}(x)=P_x \{\eta^K_hB_0^\varepsilon>0\}$, where
the $\eta^K_h$ are clusters of $\xi^K$, the truncated $K$-process of $\xi$. Fix any
bounded, uniformly continuous function $f\geq 0$ on $\RR^d$.
\vspace{-2mm}
\begin{itemize}
\item[{\em(i)}] As $0<\varepsilon^2\ll h\to
0$,
$$
\left\|\,\varepsilon^{2/\beta-d}(\beta h)^{-1/\beta}\,(p_h^\varepsilon*f)
 -c_d\,f\,\right\|\to 0. \vspace{-3mm}
$$
\item[{\em(ii)}] Fix any $b\in(0,1/2)$. Then as $0<\varepsilon^2\ll h\to
0$
with $\varepsilon^{2/\beta-d}h^{1+bd}\rightarrow 0$,
$$
\left\|\varepsilon^{2/\beta-d}\,a^{-1}_h\,(p_h^{K\varepsilon}*f)
 -c_d\,f\,\right\|\to 0. \vspace{-3mm}
$$
\end{itemize}
Both results hold uniformly over any class of uniformly bounded
and equicontinuous functions $f\geq 0$ on $\RR^d$.
\end{lemma}

{\em Proof:} (i) We follow the proof of Lemma 5.2(i) in \cite{K08}. By scaling of $\eta$ and (\ref{realone}),
\begin{eqnarray} \label{c}
\varepsilon^{2/\beta-d}(\beta h)^{-1/\beta}\lambda^d p_h^\varepsilon=(\varepsilon/\sqrt{h})^{2/\beta-d}(\beta)^{-1/\beta}\lambda^d p_1^{\varepsilon/\sqrt{h}}\rightarrow c_{\beta,d}.
\end{eqnarray}
Defining $\hat p_h^\varepsilon=p_h^\varepsilon/\lambda^dp_h^\varepsilon$, we need to show that $\|\hat p_h^\varepsilon*f-f\|\rightarrow 0$.
Writing $w_f$ for the modulus of continuity of $f$, we get
\begin{eqnarray*}
\|\hat p_h^\varepsilon*f-f\|
 &=& {\sup}_x\left|\int\hat p_h^\varepsilon(u)\,\left(f(x-u)-f(x)\right)\,du
  \right|\\
 &\leq& \int\hat p_h^\varepsilon(u)\,w_f(|u|)\,du\\
 &\leq& w_f(r)+2\,\|f\|\int_{|u|>r}\hat p_h^\varepsilon(u)\,du.
\end{eqnarray*}
It remains to show that $\int_{|u|>r}\hat p_h^\varepsilon(u)\,du\rightarrow 0$ for any fixed $r>0$. Then notice that for any fixed $r>0$ by Lemma \ref{3.hitprob}(i),
$$
\varepsilon^{2/\beta-d}(\beta h)^{-1/\beta}\int_{|u|>r} p_h^\varepsilon(u)\,du\lfrown\int_{|u|>r} p_{2h}(u)\,du\rightarrow 0.
$$

(ii) For $p^{K\varepsilon}_h$, Lemma \ref{3.hitprob}(ii) yields for any fixed $r>0$,
$$
\varepsilon^{2/\beta-d}a_h^{-1}\int_{|u|>r} p_h^{K\varepsilon}(u)\,du\lfrown\int_{|u|>r} p_{2h}(u)\,du\rightarrow 0.
$$
Following the steps of the previous proof, it is enough to show that
\begin{equation} \label{shortk}
\varepsilon^{2/\beta-d}a_h^{-1}\lambda^d p_h^{K\varepsilon}\rightarrow c_d.
\end{equation}
Since $\int_{|u|>h^b} p_{2h}(u)\,du\rightarrow 0$, Lemma \ref{3.hitprob} yields
$$
\varepsilon^{2/\beta-d}(\beta h)^{-1/\beta}1\{(B_0^{h^b})^c\}\lambda^d p_h^\varepsilon\rightarrow 0,\ \varepsilon^{2/\beta-d}a_h^{-1}1\{(B_0^{h^b})^c\}\lambda^d p_h^{K\varepsilon}\rightarrow 0.
$$
By (\ref{c}), to prove (\ref{shortk}) it suffices to show that
$$
\varepsilon^{2/\beta-d}(\beta h)^{-1/\beta}1\{B_0^{h^b}\}\lambda^d p_h^\varepsilon-\varepsilon^{2/\beta-d}a_h^{-1}1\{B_0^{h^b}\}\lambda^d p_h^{K\varepsilon}\rightarrow 0,
$$
or equivalently (by (\ref{cluster}) and (\ref{kcluster})),
$$
\varepsilon^{2/\beta-d}\left(P_{1\{B_0^{h^b}\}\lambda^d}\{\xi_h B_0^\varepsilon>0\}-P_{1\{B_0^{h^b}\}\lambda^d}\{\xi^K_h B_0^\varepsilon>0\}\right)\rightarrow 0.
$$
By Theorem 25.22 of \cite{K02} and (\ref {compensator}),
\begin{eqnarray*}
\lefteqn{\varepsilon^{2/\beta-d}\left(P_{1\{B_0^{h^b}\}\lambda^d}\{\xi_h B_0^\varepsilon>0\}-P_{1\{B_0^{h^b}\}\lambda^d}\{\xi^K_h B_0^\varepsilon>0\}\right)}\\
&\leq& \varepsilon^{2/\beta-d}E_{1\{B_0^{h^b}\}\lambda^d} N_\xi\left([0,h],(K,\infty),\RR^d)\right)\\
&=& \varepsilon^{2/\beta-d}E_{1\{B_0^{h^b}\}\lambda^d}\hat{N}_\xi\left([0,h],(K,\infty),\RR^d\right)\\
&\efrown& \varepsilon^{2/\beta-d}E\int_0^h\|\xi_s\|ds\efrown\varepsilon^{2/\beta-d}h^{1+bd}\rightarrow 0.
\end{eqnarray*}
\qed\\

\section{\large\bf Lebesgue approximations}

To prove the Lebesgue approximation for a $(2,\beta)$-process $\xi$ in $\RR^d$ with $\beta < 1$ and $d>2/\beta$, we begin with the Lebesgue approximation for $\xi^K$, the truncated $K$-process of $\xi$. Since $\xi$ and $\xi^K$ agree asymptotically as $K \rightarrow \infty$, we have thus proved the Lebesgue approximation for $\xi$. Write $\tilde c_{\beta,d}=1/c_{\beta,d}$ for convenience, where $c_{\beta,d}$ is such as in Lemma \ref{4.asymsolu}. Recall that $\xi_t^{K\varepsilon}=(\xi_t^K)^\varepsilon$, the $\varepsilon$-neighborhood measure of $\xi_t^K$.

\begin{lemma} \label{5.trunleb}
Let $\xi^K$ be the truncated $K$-process of a $(2,\beta)$-process $\xi$ in $\RR^d$ with $\beta < 1$ and $d>2/\beta$. Fix any $\mu \in \hat{\cal M}_d$ and $t>0$. Then under $P_\mu$, we have as $\varepsilon\to 0$:
$$\tilde c_{\beta,d}\,\varepsilon^{2/\beta-d}\,\xi_t^{K\varepsilon}\,\wto\,\xi^K_t \mbox{ a.s}.$$
\end{lemma}

{\em Proof:} We follow the proof of Theorem 7.1 in \cite{K08}. Fix any $f \in C_K^d$. Write $\eta_h^{Ki}$ for the
subclusters of $\xi^K_t$ of age $h$. Since the ancestors of $\xi^K_t$
at time $s=t-h$ form a Cox process directed by $\xi^K_s/a_h$, Lemma
\ref{3.variance}(i) yields
$$
E_\mu\!\left[{\sum}_i\eta_h^{Ki\varepsilon}f\,\Big|\xi^K_s\right]
 =a_h^{-1}\xi^K_s(p_h^{K\varepsilon}*f),
$$
and so by Lemma \ref{3.variance}(ii)
\begin{eqnarray*}
E_\mu\!\left|\,{\sum}_i\eta_h^{Ki\varepsilon} f
 -a_h^{-1}\,\xi_s^K(p_h^{K\varepsilon}*f)\,\right|^2
 &=& E_\mu{\rm Var}\!\left[{\sum}_i\eta_h^{Ki\varepsilon}f\,
 \Big|\xi_s^K\right]\\
 &\lfrown& a_h\,\varepsilon^{d-2/\beta}\,h^{d/2}\,\|f\|^2\,E_\mu\|\xi_s^K/a_h\|\\
 &\leq& \varepsilon^{d-2/\beta}\,h^{d/2}\,\|f\|^2\,\|\mu\|,
\end{eqnarray*}
where the last inequality follows from $E_\mu\|\xi_s^K\|\leq \|\mu\|$.
Combining with Lemma \ref{3.overlap}(i) gives
\begin{eqnarray*} \lefteqn{
E_\mu\!\left|\,\xi_t^{K\varepsilon} f
 -a_h^{-1}\,\xi_s^K(p_h^{K\varepsilon}*f)\,\right|}\qquad\\
&\leq& E_\mu\!\left|\,\xi_t^{K\varepsilon} f
 -{\sum}_i\eta_h^{Ki\varepsilon} f\,\right|
 +E_\mu\!\left|\,{\sum}_i\eta_h^{Ki\varepsilon} f
 -a_h^{-1}\,\xi_s^K(p_h^{K\varepsilon}*f)\,\right|\\
&\lfrown& \varepsilon^{2(d-2/\beta)}\,h^{1-d/2}\,\|f\|
 +\varepsilon^{1/2(d-2/\beta)}\,h^{d/4}\,\|f\|\\
&=& \varepsilon^{d-2/\beta}\left(\varepsilon^{d-2/\beta}h^{1-d/2}+\varepsilon^{-1/2(d-2/\beta)}h^{d/4}
 \right)\|f\|.
\end{eqnarray*}
Let $c$ satisfy
\begin{equation} \label{definition}
(d-2/\beta)+(-d/2+1/2)c=0.
\end{equation}
Clearly $c\in(0,2)$. Taking $\varepsilon=r^n$ for a fixed $r\in(0,1)$ and $h=\varepsilon^c=r^{cn}$, and writing
$s_n=t-h=t-r^{cn}$, we obtain
\begin{eqnarray*} \lefteqn{
E_\mu\,\sum_nr^{n(2/\beta-d)}\left|\,\xi_t^{Kr^n}\!f-
 a_{r^{cn}}^{-1}\,\xi^K_{s_n}(p_{r^{cn}}^{Kr^n}*f)\,\right|}\qquad\\
&& \lfrown\;
 \sum_n\!\left(r^{[(d-2/\beta)+(-d/2+1)c]n}+r^{[-1/2(d-2/\beta)+(d/4)c]n}\right)\|f\|<\infty,
\end{eqnarray*}
since $(d-2/\beta)+(-d/2+1)c>0$ and $-1/2(d-2/\beta)+(d/4)c>0$ by (\ref {definition}). This implies
\begin{equation} \label{thefirstterm}
r^{n(2/\beta-d)}\left|\,\xi_t^{Kr^n}\!f
 -a_{r^{cn}}^{-1}\,\xi^K_{s_n}(p_{r^{cn}}^{Kr^n}*f)\,\right|
 \to 0\,\mbox{ a.s. }P_\mu.
\end{equation}

Now we write
\begin{eqnarray*}
\lefteqn{\left|\,\varepsilon^{2/\beta-d}\,\xi_t^{K\varepsilon} f-c_{\beta,d}\,\xi^K_tf\,
  \right|}\\
  &\leq& \varepsilon^{2/\beta-d}\left|\,\xi_t^{K\varepsilon} f-a_h^{-1}\,\xi^K_s(p_h^{K\varepsilon}*f)\,\right|+c_{\beta,d}\,|\xi^K_sf-\xi^K_tf|\\
&&+\;\|\xi^K_s\|\left\|\,\varepsilon^{2/\beta-d}\,a_h^{-1}\,(p_h^{K\varepsilon}*f)
  -c_{\beta,d} f\,\right\|.
\end{eqnarray*}
For the last term, we first fix $b=1/2-1/d$, then apply Lemma \ref{4.shorttime}. Noting that by (\ref {definition})
$$(2/\beta-d)+(1+bd)c=(2/\beta-d)+(d/2)c>0,$$
we get by Lemma \ref{4.shorttime}
$$\left\|\,\varepsilon^{2/\beta-d}\,a_h^{-1}\,(p_h^{K\varepsilon}*f)-c_{\beta,d} f\,\right\|\rightarrow 0$$
along the sequence $(r^n)$. Using (\ref {thefirstterm}) and the
a.s.\ weak continuity of $\xi^K$ at the fixed time $t$, we see that the right-hand side tends a.s.\ to 0 as
$n\to\infty$, which implies $\varepsilon^{2/\beta-d}\,\xi_t^{K\varepsilon} f\rightarrow c_{\beta,d}\,\xi^K_tf$
 a.s.\ as $\varepsilon\to 0$ along the sequence
$(r^n)$ for any fixed $r\in(0,1)$. Since this holds
simultaneously, outside a fixed null set, for all rational
$r\in(0,1)$, the a.s.\ convergence extends by Lemma 2.3 in \cite{K08}
to the entire interval $(0,1)$.

Applying this result
to a countable, convergence-determining class of functions $f$
(cf.\ Lemma 3.2.1 in \cite{D93}), we obtain the required a.s.\
vague convergence. Since $\mu$ is finite, the $(2,\beta)$-process $\xi_t$ has a.s.\
compact support (cf. Theorem 9.3.2.2 of \cite{D93} and the proof of Theorem 1.2 in \cite{DIP89}). By Lemma \ref{2.domination}, $\xi^K_t$ also has a.s.\
compact support, and so the a.s.
convergence remains valid in the weak sense.
\qed\\

Now we may prove our main result, the Lebesgue approximation of $(2,\beta)$-processes. Again, we write $\tilde c_{\beta,d}=1/c_{\beta,d}$ for convenience, where $c_{\beta,d}$ is such as in Lemma \ref{4.asymsolu}. Also recall that $\xi_t^\varepsilon=(\xi_t)^\varepsilon$ denotes the $\varepsilon$-neighborhood measure of $\xi_t$. For random measures $\xi_n$ and $\xi$ on $\RR^d$, $\xi_n\vto\xi$ (or $\wto$) in $L^1$ means that $\xi_n f\to\xi f$ in $L^1$ for all $f$ in $C^d_K$ (or $C^d_b$).

\begin{theorem} \label{5.lebesgue}
Let the $(2,\beta)$-process $\xi$ in $\RR^d$ with $\beta<1$ and $d>2/\beta$ be locally finite under
$P_\mu$, and fix any $t>0$. Then under $P_\mu$, we have as
$\varepsilon\to 0$:
$$
\tilde c_{\beta,d}\,\varepsilon^{2/\beta-d}\,\xi_t^\varepsilon\vto \xi_t \mbox{ a.s and in }L^1.
$$
This remains true in the weak sense when $\mu$ is finite. The weak version holds even for the clusters $\eta_t$ when $\|\mu\|=1$.
\end{theorem}

{\em Proof:}
For a finite initial measure $\mu$, by Lemma \ref{5.trunleb} and Lemma \ref{2.asymsame} we get as $\varepsilon\to 0$
$$
\tilde c_{\beta,d}\,\varepsilon^{2/\beta-d}\,\xi_t^\varepsilon\wto \xi_t \mbox{ a.s}.
$$
For a general $\mu\in{\cal M}_d$ with $\mu p_t<\infty$ for
all $t>0$, write $\mu=\mu'+\mu''$ for a finite $\mu'$, and let
$\xi=\xi'+\xi''$ be the corresponding decomposition of $\xi$ into
independent components with initial measures $\mu'$ and $\mu''$.
Fixing an $r>1$ with ${\rm supp}\,f\subset B_0^{r-1}$ and using
the result for finite $\mu$, we get a.s.\ on $\{\xi''_tB_0^r
=0\}$
$$
\varepsilon^{2/\beta-d}\,\xi_t^\varepsilon f
 =\varepsilon^{2/\beta-d}\,\xi_t'^\varepsilon f
 \to c_{\beta,d}\,\xi'_tf=c_{\beta,d}\,\xi_tf.
$$
As $\mu'\uparrow\mu$, we get by Lemma \ref{3.initial}
$$
P_\mu\{\xi''_tB_0^r=0\}=P_{\mu''}\{\xi_tB_0^r=0\}\to 1,
$$
and the a.s.\ convergence extends to $\mu$. As in the proof of Lemma \ref{5.trunleb},
we can obtain the required a.s. vague convergence.

To prove the convergence in $L^1$, we note that for any $f\in
C_K^d$
\begin{eqnarray}\label{l1convergence}
\varepsilon^{2/\beta-d}E_\mu\xi_t^\varepsilon f
 &=& \varepsilon^{2/\beta-d}\int P_\mu\{\xi_tB_x^\varepsilon>0\}\,f(x)\,dx
 \nonumber\\
 &\to& \int c_{\beta,d}\,(\mu*p_t)(x)\,f(x)\,dx
 =c_{\beta,d}\,E_\mu\xi_tf,
\end{eqnarray}
by Theorem \ref{4.asymhit}. Combining this with the a.s.\
convergence under $P_\mu$ and using Proposition 4.12 in
\cite{K02}, we obtain $E_\mu |\varepsilon^{2/\beta-d}
\,\xi_t^\varepsilon f-c_{\beta,d}\,\xi_tf|\to 0$. For finite $\mu$,
(\ref{l1convergence}) extends to any $f\in C_b^d$ by dominated
convergence based on Lemmas \ref{3.hitcomp} and \ref{3.hitprob}(i),
together with the fact that $\lambda^d (\mu*p_t)=\|\mu\|<\infty$
by Fubini's theorem.

To extend the Lebesgue approximation to the individual clusters $\eta_t$, let
$\zeta_0$ denote the process of ancestors of $\xi_t$ at time 0,
and note that
$$
P_x\{\eta_t\in\cdot\}=P_{\delta_x}[\xi_t\in\cdot|\|\zeta_0\|=1],
$$
where $P_{\delta_x}\{\|\zeta_0\|=1\}=(\beta t)^{-1/\beta}e^{-(\beta t)^{-1/\beta}}>0$. The a.s.
convergence then follows from the corresponding statement for
$\xi_t$. Since
$$P_\mu\{\eta_t\in\cdot\}=\int\mu(dx)P_x\{\eta_t\in\cdot\},$$
the a.s. convergence under any $P_\mu$ with $\|\mu\|=1$ also follows.
To obtain the weak $L^1$-convergence in this case, we note that for $f\in C_b^d$,
\begin{eqnarray*}
\varepsilon^{2/\beta-d}E_\mu\eta_t^\varepsilon f
 &=& \varepsilon^{2/\beta-d}\int P_\mu\{\eta_tB_x^\varepsilon>0\}\,f(x)\,dx
 \nonumber\\
 &\to& c_{\beta,d}\,(\beta t)^{1/\beta}\int(\mu*p_t)(x)\,f(x)\,dx
 =c_{\beta,d}\,E_\mu\eta_t f,
\end{eqnarray*}
by dominated convergence based on Lemma \ref{3.hitprob}(i) and Theorem \ref{4.asymhit}.
\qed\\

As in Corollary 7.2 of \cite{K08}, for the intensity measures in Theorem \ref{5.lebesgue}, we have even convergence in total variation.

\begin{corollary}
Let $\xi$ be a $(2,\beta)$-process in $\RR^d$ with $\beta<1$ and $d>2/\beta$. Then for any finite $\mu$ and $t>0$, we have as $\varepsilon\to 0$$\rm:$
$$
\left\|\,\varepsilon^{2/\beta-d}E_\mu\xi_t^\varepsilon
 -c_{\beta,d}\,E_\mu\xi_t\,\right\|\to 0.
 $$
This remains true for the clusters $\eta_t$, and it also
holds locally for $\xi_t$ whenever $\xi$ is locally finite under
$P_\mu$.
\end{corollary}

{\em Acknowledgement:}\/ This paper is based on part of my Ph.D. dissertation at Auburn University. I sincerely thank my advisor Olav Kallenberg for his encouragement and suggestions.\\

{\small

\end{document}